\documentclass[12pt]{article}
\usepackage[left=2.5cm,right=2.5cm,top=3cm,bottom=4cm]{geometry}
\usepackage{authblk} % Este paquete permite escribir varios autores.
\usepackage{verbatim}
\usepackage{framed}
\usepackage{booktabs}%para toprule
\usepackage{amsmath,amssymb,amsfonts,latexsym,cancel}
\usepackage{amsthm}
\usepackage{relsize,exscale}
\usepackage{enumerate}
\usepackage{color}
\usepackage{xcolor}
\usepackage{graphics}
\usepackage{float}
\usepackage{array}
\usepackage{etoolbox}
\usepackage{amsbsy}
\usepackage{tikz}
\usetikzlibrary{babel}
\usetikzlibrary{arrows,chains,matrix,positioning,scopes}

\makeatletter
\tikzset{join/.code=\tikzset{after node path={%
\ifx\tikzchainprevious\pgfutil@empty\else(\tikzchainprevious)%
edge[every join]#1(\tikzchaincurrent)\fi}}}
\makeatother
\tikzset{>=stealth',every on chain/.append style={join},
         every join/.style={->}}
\tikzstyle{labeled}=[execute at begin node=$\scriptstyle,
   execute at end node=$]
\newtheorem{theorem}{Theorem}[section]

\newtheorem{proposition}[theorem]{Proposition}
\newtheorem{remark}[theorem]{Remark}
\newtheorem{corollary}[theorem]{Corollary}
\newtheorem{definition}[theorem]{Definition}

\title{Approximations of the Ruin Probability\\
in a Discrete Time Risk Model}
\author[1]{David J. Santana}
\author[2]{Luis Rinc\' on}
\affil[1]{Divisi\'on Acad\'emica de Ciencias B\'asicas\\UJAT\\M\'exico}
\affil[2]{Departamento de Matem\'aticas\\Facultad de Ciencias\\UNAM\\M\'exico}

\begin{document}
\maketitle

\begin{abstract}
Based on a discrete version of the Pollaczeck-Khinchine formula, a general
method to calculate the ultimate ruin probability in the Gerber-Dickson
risk model is provided when claims follow a negative binomial mixture
distribution. The result is then extended for claims with a mixed Poisson
distribution. The formula obtained allows for some approximation procedures.
Several examples are provided along with the numerical evidence of the
accuracy of the approximations.
\end{abstract}

%-----------------------------------------------------------------
\section{Introduction}
\begin{figure}
\centering
\input{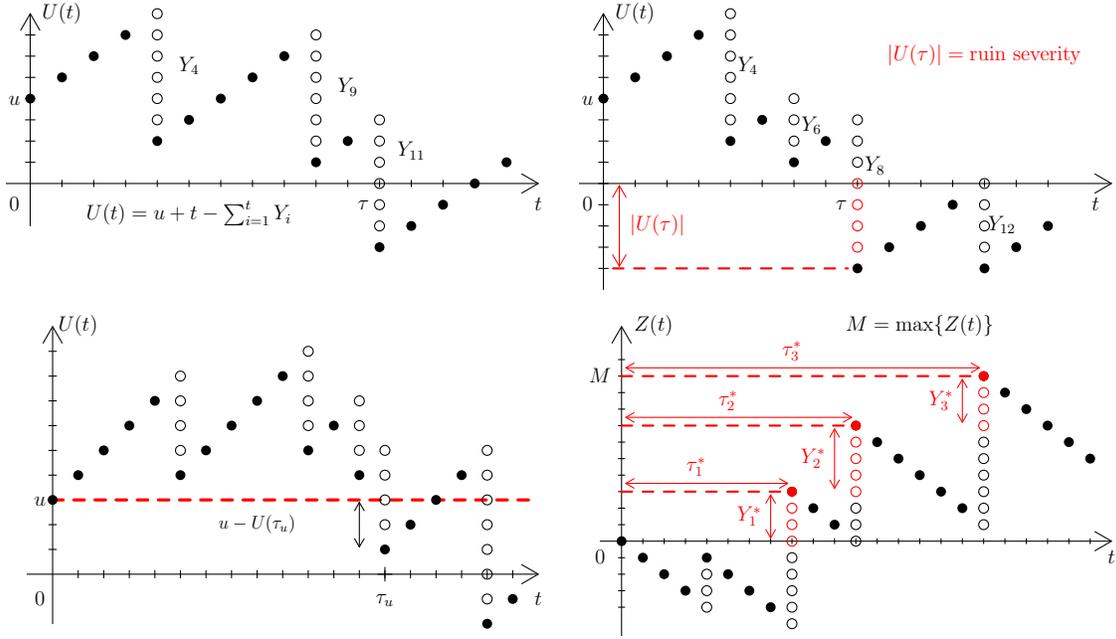}
\caption{Discrete time risk process trajectories and some related quantities.} \label{plots1}
\end{figure}

Several models have been proposed
%versions
for a discrete time\footnote{We will reserve the use of letter $n$ for the
approximation procedures proposed later on.}
risk process $\{U(t):t=0,1,\ldots\}$.
%The use of the
%letter $ n $ has been avoided to denoting discrete time, as traditionally used, because this letter will be used in approximations, both to the mixed Poisson distribution functions and to the ruin probability. This convention is followed throughout the work.}
%, depending on the assumptions made for the premiums, the amount of the claims and the number of claims per period. In 
The following model is known as a compound binomial process and
was first considered in~\cite{1988gerber},
\begin{equation}
\label{Compound-binomial-process}
U (t)=u + t - \sum_{i = 1}^{N (t)} X_i ,
\end{equation}
where $U(0)=u\ge0$ is an integer representing the initial capital
and the counting
process $\{N(t):t=0,1,\ldots\}$ has a $\mbox{Binomial}(t,p)$ distribution,
%being the number of periods where claims have been made up to time $t$ and has distribution $\mbox{bin}(t,p) $,
where $p$ stands for the probability of a claim in each period.
The discrete random variables $X_1, X_2, \ldots $ are i.i.d. with
probability function $ f_X(x) = P (X_i = x) $ for $x= 1,2,\ldots $ and
mean $\mu_X$ such that $\mu_X\cdot p<1$.
This restriction comes from the net profit condition.
Each $X_i $ represents the total amount of claims in the $i$-th period where
claims existed. In each period, one unit of currency from premiums is gained.
The top-left plot of Figure~\ref{plots1} shows a realization of this risk process.
%This model .
%assumes that all elements are nonnegative integers, so that the risk process 
%$\{U(t):t=0,1,\ldots\}$ has discrete state space.
The ultimate ruin time is defined as
\begin{equation*}
\tau=\min\,\{t \ge 1: U (t) \le 0 \},
\end{equation*}
as long as the indicated set is not empty, otherwise $\tau:=\infty$.
Hence, the probability of ultimate ruin is
\begin{equation*}
\psi (u)=P (\tau <\infty \mid U (0) = u).
\end{equation*}
One central problem in the theory of ruin is to find $\psi (u)$. For the above
model this probability can be calculated using the following
relation known as Gerber's formula~\cite{1988gerber},
\begin{eqnarray}
\psi (0) &=& p\cdot\mu_X, \\
\label{gerber-ruin}
\psi (u) &=& (1-p) \psi (u + 1) + p \sum_{x = 1}^u \psi (u + 1-x) \, f_X(x) + p \, \overline{F}_X(u),
\end{eqnarray}
for $u=1,2,\ldots$ where
$\overline{F}_X(u)=P(X_i>u)=\sum_{x= u + 1}^\infty f_X(x)$. 

*****

% ----------------------------------------
An apparently simpler risk model is defined as follows.

\begin{definition}
Let $u\ge0$ be an integer and let $Y_1, Y_2,\ldots$ be i.i.d. random variables
taking values in $\{0,1,\ldots\}$. The Gerber-Dickson risk process
$\{U (t):t=0,1,\ldots \}$ is given by
\begin{equation}
\label{defmd}
U (t) = u + t- \sum_{i = 1}^t Y_i.
\end{equation}
\end{definition}
%The ruin time $ \tau $ and the ruin probability with infinite horizon are defined as before.

%\begin{equation}\label{57}
%\tau = \min \{t \ge 1: U (t) \le 0 \},
%\end{equation}
%assuming that the indicated set is different from the empty, otherwise $ \tau: = \infty $.

%\begin{equation}
%\psi (u) = P (\tau <\infty \mid U (0) = u),
%\end{equation}
%for each value of $ u = 0,1,2, \ldots $ For the net profit condition to be fulfilled, it is required that,

In this case, at each unit of time there is always a claim of size $Y$.
If $\mu_Y$ denotes the expectation of this claim, the net profit condition now
reads $\mu_Y<1$.  It can be shown~\cite[ pp. 467] {2015damarackas} that this
condition implies that $\psi(u)<1$, where the time of ruin $\tau$ and the
ultimate ruin probability $\psi(u)$ are defined as before.
Under a  conditioning argument
it is easy to show that the probability of ruin satisfies the recursive relation

\begin{eqnarray}
\label{psicerod}
\psi(0)&=&\mu_Y,\\
\label{psid}
\psi(u)&=&\sum\limits_{y=0}^u {f_Y(y)\,\psi (u + 1-y)} + \overline{F}_Y(u),
\quad u\ge1.
\end{eqnarray}

% ------------------------------------------------------------
Now, given a compound binomial model~(\ref{Compound-binomial-process})
we can construct a Gerber-Dickson model~(\ref{defmd}) as follows.
Let $R_1,R_2,\ldots$ be i.i.d. $\mbox{Bernoulli}(p)$ random variables and
define $Y_i=R_i\cdot X_i$, $i\ge1$. The distribution of these claims is
$f_Y(0)=1-p$ and $f_Y(y)=p\cdot f_X(y)$ for $y\ge1$.

Conversely, given model~(\ref{defmd}) and defining $p=1-f_Y(0)$, we
can construct a
model~(\ref{Compound-binomial-process}) by letting claims $X_i$ have
distribution $f_X(x)=f_Y(x)/p$, for $x\ge1$. It can be readily checked that
$\mu_Y=p\cdot\mu_X$ and that the probability generating function
of $U(t)$ in both models coincide. This shows
models~(\ref{Compound-binomial-process}) and~(\ref{defmd}) are
equivalent in the sense that $U(t)$ has the same distribution in both models.
As expected, the recursive relations~(\ref{gerber-ruin})
and~(\ref{psid}) can be easily obtained one from the other.

In this work we will use the notation in the Gerber-Dickson risk
model~(\ref{defmd}) and drop the subindex in the distribution of claims.
Also, as time an other auxiliary variables are considered discrete, we will
write, for example, $t\ge0$ instead of $t=0,1,\ldots$
Our main objective es to provide some methods to approximate the
ultimate ruin probability in the discrete risk model of Gerber and Dickson.

A survey of results and models for discrete time risk models can be found
in~\cite{2009Li-Lu-Garrido}.

%-----------------------------------------------------------------

\section{The Pollaczeck--Khinchine formula}
\label{five}

The continuous version of this formula plays a major role in the theory of ruin
for the Cram\'er-Lundberg model. On the contrary, the discrete version is
seldom mentioned in the literature on discrete time risk models.
In this section we develop this formula and apply it later to find a general
method to calculate ultimate ruin probabilities for claims with particular distributions.
The construction procedure resembles closely that for the continuous case.

Assuming $\tau<\infty$, the non-negative
random variable $W=|U(\tau)| $ is known as the severity of ruin.
It indicates how large the capital drops below zero at the time of ruin.
See the top-right plot of Figure~\ref{plots1}.
%Note that there can be zero severity due to the way of defining the ruin time in (\ref{57}). In \cite{1988gerber} some results related to the distribution of $Y$ are obtained.
\begin{comment} ***
\\ \textcolor{red}{***}\\
In~\cite{2000cheng}, \cite{1994dickson}, \cite{2002Li} and \cite{2005li} the Gerber-Shiu function is studied in its discrete version and some results are
obtained regarding the joint distributions, the marginal distributions and the moments of $ Y $, of the capital just before the ruin $ | U (\tau - 1) | $ and of the ruin time.
\\ \textcolor{red}{***}\\
\end{comment}
The joint probability of ruin and severity not greater than $w= 0,1,\ldots $
is denoted by
\begin{equation}\label{defin}
\varphi(u,w) = P (\tau <\infty, W\le w \mid U (0) = u).
\end{equation}
In~\cite{1994dickson} it is shown that, in particular,
\begin{equation}\label{propos}
\varphi (0, w) = \sum_{x = 0}^{w}\overline{F}(x),\quad w\ge0.
\end{equation}
Hence,
\begin{equation}
\label{seve-proba}
P (\tau<\infty, W=w\mid U (0) = 0)=\varphi (0,w)-\varphi (0,w-1)
=\overline{F}(w).
\end{equation}

This probability will be useful in finding the distribution
of the size of the first drop of the risk process below its initial capital $u$,
see Proposition \ref{fdeLu} below, which will ultimately lead us to the Pollaczeck--Khinchine formula.
For every claim distribution, there is an associated distribution which often
appears in the calculation of ruin probabilities. This is defined next.

\begin{definition}
Let $ F (y) $ be the distribution function of a discrete random variable
with values ​​$0,1,\ldots $ and with finite mean $\mu\ne0$. Its equilibrium probability function is defined by
\begin{equation}
\label{515}
f_e (y) = \overline{F}(y)/\mu,\quad y\ge0.
\end{equation}
\end{definition}
The probability function defined by~(\ref{515}) is also known as the
integrated-tail distribution, although this name is best suited to
continuous distributions.
For example, the equilibrium distribution associated to a $\mbox{Geometric}(p)$
claim distribution with mean $\mu=1/(1-p)$ is the same geometric since
\begin{equation}
\label{532}
f_e (y) = \overline{F}(y)/\mu=(1-p)^{y +1} \, p / (1-p) = p \, (1-p)^{y},
\quad y\ge0.
\end{equation}

As in the continuous time risk models, let us define
the surplus process $\{Z (t):t\ge0\}$ by
\begin{equation}
\label{procesoZn}
Z (t)=u-U (t)=\sum_{i = 1}^{t}(Y_i-1).
\end{equation}
This is a random walk that starts at zero, it has stationary and independent
increments and $Z (t)\rightarrow-\infty\ a.s.$ as $t\to\infty$ under the
net profit condition $\mu<1$. See bottom-right plot of Figure~\ref{plots1}.
In terms of this surplus process,
ruin occurs when $Z(t)$ reaches level $u$ or above. Thus, the ruin probability
can be written as
\begin{equation}
\label{eye-d-0}
\psi (u) = P (Z (t) \ge u \mbox{ for some } t \ge 1) = P \left(\max \limits_{t \ge 1} \left\{Z (t) \right\} \ge u \right),\quad u\ge1.
\end{equation}

As $u\ge1$ and $Z(0)=0$, we can also write
\begin{equation}
\label{eye-d}
\psi (u) = P \left(\max \limits_{t \ge 0} \left\{Z (t) \right\} \ge u \right).
\end{equation}

We next define the times of records and the severities for the surplus process.

\begin{definition}
Let $\tau_0^* := 0$. For $i\ge1$ the $i$-th record time of the
surplus process is defined as
\begin{equation}
\label{records} 
\tau_i^* = \min\,\{t> \tau_{i-1}^*: Z (t) \ge Z(\tau_{i-1}^*) \},
\end{equation}
when the indicated set is not empty, otherwise $\tau_i^* :=\infty$.
The non-negative variable $Y_i^*=Z(\tau_i^*)-Z (\tau_{i-1}^*)$ is called the
severity or size of the $i$-th record time, assuming $\tau_i^*<\infty$.
\end{definition}

The random variables $\tau_0^*<\tau_1^*<\cdots$ represent the
stopping times when the surplus process
$\left\{Z (t):t\ge0 \right\}$ arrives at a new or the previous
maximum, and the severity $Y_i^*$ is the difference between the maxima at
$\tau_i^*$ and $\tau_{i-1}^*$.
A graphical example of these record times are shown in the bottom-right plot of
Figure~\ref{plots1}.
%\begin{remark}\label{522}
%If $\tau_k^* <\infty $ for some $ k $, then $ \tau_i^* <\infty $ for $ i = 1,2, \ldots, k-1 $, and if $ \tau_m^* = \infty $, then $ \tau_i^* = \infty $ for $ i = m + 1, m + 2, \ldots $ These random variables are stop times with respect to the natural filtration of the process $ \left\{Z (t) \right\}_{t = 0}^\infty $, that is, its value at $t$ can be determined by the trajectory of the process up to $t$.
%\end{remark}
In particular, observe $\tau_1^*$ is the first positive time the risk process is
less than or equal to its initial capital $u$, that is,
\begin{equation}
\label{fu}
\tau_1^*=\min\,\{t>0:u-U (t)\ge0\},
\end{equation}
and the severity is $Y_1^*=Z(\tau_1^*)=u-U(\tau_1^*)$ and this is the size
of this first drop below level $u$.
Also, since the surplus process has stationary
increments, all severities share the same distribution, that is,
\begin{equation}
\label{severities-share-distribution}
Y_i^*=Z (\tau_i^*) - Z (\tau_{i-1}^*)\sim Z (\tau_1^*) - Z (0)=Y_1^*,\quad i\ge1,
\end{equation}
assuming $\tau_i^*<\infty$. We will next find out that distribution.

% \textcolor{red}{Observemos que $\tau_u$ es un tiempo de ruina pero
% cuando la ruina se redefine como $(U (t)\le u)$.}

\begin{comment} *****
Comparing the definitions of the time of the first fall and the time of ruin, it is clear that the time of the first fall is equal to the time of ruin when the initial capital is zero. Using the respective notation,
\begin{equation}
P (\tau_0 = t \mid \tau_0 <\infty) = P (\tau = t \mid \tau <\infty, U (0) = 0).
\end{equation}

An important result used later is that the probability function of the random variable $ Lu $ is equal to the equilibrium probability function of claims. The following proposition states this fact and its demonstration is a simple application of conditional probability.
\end{comment} *****

\begin{proposition}
\label{fdeLu}
Let $k\ge1$.
Conditioned on the event $(\tau_k^*<\infty)$, the severities $Y_1^*,\ldots,Y_k^*$ are
independent and identically distributed according to the equilibrium distribution
\begin{equation}
\label{feqd}
P(Y^*=x \mid \tau_k^*<\infty)=\overline{F}(x)/\mu,
\quad x\ge0.
\end{equation}
\end{proposition}

\textbf{Proof.}
By~(\ref{severities-share-distribution}), it is enough to find the distribution of $Y_1^*$.
Observe that $\tau_1^*=\tau$ when $U(0)=0$.
By~(\ref{seve-proba}) and~(\ref{dickson-ruin-0}), for $x\ge0$,
\begin{eqnarray*}
P (Y_1^*=x \mid \tau_1^*<\infty))
&=& P (u - U (\tau_1^*) = x \mid \tau_1^*<\infty) \\
%&=& P (-U (\tau_0) = y \mid \tau_0 <\infty) \\
&=& P ( | U (\tau) | = x \mid \tau <\infty, U (0) = 0) \\
&=& P (\tau <\infty, Y = x \mid U (0) = 0) / P (\tau <\infty \mid U (0) = 0) \\
&=& \overline{F}(x) / \mu.
\end{eqnarray*}
The independence property follows from the independence of the claims. Indeed,
the severity of the $i$-th record time is
\begin{equation*}
Y_i^*=Z (\tau_i^*) - Z (\tau_{i-1}^*) = \sum_{j = \tau_{i-1}^* + 1}^{\tau_{i}^*}( Y_j -1),\quad i\ge1.
\end{equation*} 
Therefore,
\begin{displaymath}
P \left(\bigcap_{i = 1}^k \left( Y_i^* = y_i \right)\right) 
=P \left(\bigcap_{i = 1}^k \left(\sum_{j = \tau_{i-1}^* + 1}^{\tau_{i}^*}( Y_j -1) = y_i \right) \right) 
=\prod_{i = 1}^k \, P \left( Y_i^* = y_i \right).
\end{displaymath}
\begin{comment} ***
We first consider the case $i=1$. Observe that $\tau_1^*=\tau$ when $U(0)=0$.
By~(\ref{seve-proba}) and~(\ref{dickson-ruin-0}), for $x\ge0$,
\begin{eqnarray*}
P (Y_1^*=x \mid \tau_1^*<\infty))
&=& P (u - U (\tau_1^*) = x \mid \tau_1^*<\infty) \\
%&=& P (-U (\tau_0) = y \mid \tau_0 <\infty) \\
&=& P ( | U (\tau) | = x \mid \tau <\infty, U (0) = 0) \\
&=& P (\tau <\infty, Y = x \mid U (0) = 0) / P (\tau <\infty \mid U (0) = 0) \\
&=& \overline{F}(x) / \mu.
\end{eqnarray*}
%The last equality occurs thanks to (\ref{psicerod}) and (\ref{seve-proba}). 
For $i\ge2$ define the event $A_i(s, y)=(\tau_{i-1}^* = s, Z (\tau_{i-1}^*) = y)$,
for $s\ge i-1 $ and $y\ge0$. Then
\begin{eqnarray}
\nonumber
P(Y_i^*=x \mid \tau_i^*<\infty)
&=& \sum_{s = i-1}^\infty \sum_{y = 0}^\infty P (Z (\tau_i^*) - Z (\tau_{i-1}^*) = x \mid A_i(s, y), \tau_i^*<\infty) \, P(A_i(s, y)) \\
\nonumber &=& \sum_{s = i-1}^\infty \sum_{y = 0}^\infty P (Z (\tau_i^*) - Z (s) = x \mid A_i(s, y), \tau_i^*<\infty) \, P(A_i(s, y)) \\
\nonumber &=& \sum_{s = 0}^\infty \sum_{y = 0}^\infty P (Z (\tau_1^*) = x) \, P(A_i(s, y)) \\
\label{stationary-increments-of-Ztaui-d}
&=& \overline{F}(x)/\mu.
\end{eqnarray}
%The third equality of the previous development occurs thanks to the stationary and independent increases of the surplus process.
\end{comment}
\qed

\begin{comment}
If $ U (0) = 0 $, then the definitions of $ \tau_1^* $, $ \tau_0 $ and the time of ruin $ \tau $ coincide.
\begin{equation}\label{observation215-d}
\tau_1^* = \min \{t> 0: Z (t) \ge 0 \} = \min \{t \ge 1: U (t) \le U (0) = u \} = \tau_u.
\end{equation}
Therefore, $ Z (\tau_1^*) = Z (\tau_u) = u-U (\tau_u) $, then, $Z (\tau_1^*) = Lu$, and by (\ref{feqd}), the probability function of $ Z (\tau_1^*) $ is equal to the equilibrium probability function $ f_e $ defined in (\ref{515}).
\end{comment}

Since the surplus process is a Markov process, the following properties hold:
For $i\ge2$ and assuming $\tau_i^*<\infty $, for $ 0<s<x $,
\begin{equation}
\label{observation-2-tau-star-d}
P (\tau_i^* = x \mid \tau_{i-1}^* = s) = P (\tau_i^* - \tau_{i-1}^* = x - s \mid \tau_{i-1}^* = s)  = P (\tau_1^* = x - s).
\end{equation}
Also, for $k\ge1$,
\begin{eqnarray}
\label{observation-2-tau-star-d-1}
P (\tau_k^* <\infty \mid \tau_{k-1}^{*} <\infty)&=&P (\tau_1^* <\infty),\\
\label{observation-2-tau-star-d-2}
P (\tau_k^* =\infty \mid \tau_{k-1}^{*} <\infty)&=&P (\tau_1^* =\infty).
\end{eqnarray}

The total number of records of the surplus process $\{Z (t):t\ge0\}$ is defined by
the non-negative random variable
\begin{equation}
\label{number-of-records-d}
K = \max\,\{k \ge 1: \tau_k^* <\infty \},
\end{equation} 
when the indicated set is not empty, otherwise $K:=0$.
Note that $0\le K<\infty$ a.s. since $ Z (t) \rightarrow - \infty $ a.s. under the net profit condition. The distribution of this random variable is established next.

\begin{proposition}
\label{num-de-jumps-d}
The number of records $K$ has a $\mbox{Geometric}(1-\mu)$ distribution,
that is,
\begin{equation}
\label{tests-of-kd}
f_K (k)=(1-\mu)\mu^k,\quad k\ge0.
\end{equation}
\end{proposition}

\textbf{Proof.} The case $k=0$ can be related to the ruin probability with $u=0$
as follows,
\begin{equation*}
f_K (0) = P (\tau_1^{*}=\infty)
=P(\tau = \infty \mid U(0) = 0)=1-\psi(0)=1-\mu.
\end{equation*}
Hence, $P(K>0)=\psi (0)=\mu$. Let us see the case $k=1$,
\begin{equation*}
f_K (1) = P (\tau_1^{*} <\infty, \tau_2^{*} = \infty) 
= P ( \tau_2^{*} = \infty \mid \tau_1^{*} <\infty) P( \tau_1^{*} <\infty).
\end{equation*}
By~(\ref{observation-2-tau-star-d-1}),
\begin{equation*}
f_K (1) = P (\tau_1^{*} = \infty) P (\tau_1^{*} <\infty) 
= P (K> 0) f_K(0)
=\mu(1-\mu).
\end{equation*}
Now consider the case $k\ge2$ and let $A_k=(\tau_k^* <\infty)$.
Conditioning on $A_{k-1}$ and its complement,
\begin{eqnarray*}
P(A_k) &=& P (\tau_k^* <\infty \mid A_{k-1}) P ( A_{k-1}) \\
&=& P (\tau_k^* <\infty \mid \tau_{k-1}^{*} <\infty) P ( A_{k-1}) \\
&=& P (\tau_1^* <\infty) P ( A_{k-1}) \\
&=& \psi (0) P ( A_{k-1}).
\end{eqnarray*}

An iterative argument shows that $P(A_k)=(\psi(0))^k$, $k\ge2$.
Therefore,
\begin{equation*}
f_K (k) = P (\tau_{k + 1}^{*} = \infty, A_k ) 
= P (\tau_{k + 1}^{*} = \infty \mid A_k ) P ( A_k ) 
= P (\tau_1^{*} = \infty) (\psi (0))^k 
= (1-\mu)\mu^k.
\end{equation*}
\qed

In the following proposition it is established that the ultimate maximum
of the surplus process
%$\max\limits_{t \ge 1}\,\{Z (t) \}$
has a compound geometric distribution. This will allow us to write the ruin
probability as the tail of this distribution.
%, with $K$ as the counting variable
%and summands the variables $Y_i^*$, $i\ge1$.
\begin{proposition}
\label{kym}
For a surplus process $\{Z (t):t\ge0\}$ with total number of records $K\ge0$ and record severities
$ Y_1^*, Y_2^*, \ldots, Y_K^* $, 
\begin{equation}
\max \limits_{t\ge0}\,\{Z (t) \}\stackrel{d}{=}  \sum_{i=1}^{K} Y_i^*.
\end{equation}
Hence,
\begin{equation}
\label{psi-geo-comp-d}
\psi (u) = P \left(\sum_{i = 1}^{K} Y_i^* \ge u \right),\quad u\ge1.
\end{equation}
\end{proposition}
\textbf{Proof.}
\begin{equation}
\label{eye2-d}
\sum_{i =1}^{K} Y_i^*=\sum_{i = 1}^{K} \left(Z (\tau_i^*) - Z (\tau_{i-1}^* ) \right)
=Z (\tau_K^*)=\max_{t\ge0}\,\{Z (t)\}\quad a.s.
\end{equation}
Thus, for $u\ge1$,
\begin{equation*}
\psi(u) = P \left(\max_{t\ge0}\,\{Z (t)\}\ge u \right) 
= P \left(\sum_{i = 1}^{K} Y_i^*\ge u \right).
\end{equation*}
\qed

\begin{proposition}
\label{PKd}
(Pollaczeck--Khinchine formula, discrete version)
The probability of ruin for a Gerber-Dickson risk process can be written as
\begin{equation}
\label{PK}
\psi (u) = (1-\mu) \sum_{k = 1}^{\infty} P (S_k^* \ge u) \, \mu^k,\quad u\ge0,
\end{equation}
where $ S_k^* = \sum_{i = 1}^{k} Y^*_i$.
\end{proposition}
\textbf{Proof.}
For $u=0$, the sum in~(\ref{PK}) reduces to $\mu$ which we know is $\psi(0)$.
For $u\ge1$, by~(\ref{tests-of-kd}) and~(\ref{psi-geo-comp-d}),
\begin{equation*}
\psi (u) = P \left(\sum_{i = 1}^{K} Y^*_i \ge u \right) 
= \sum_{k = 0}^{\infty} P \left(\sum_{i = 1}^{K} Y^*_i \ge u \mid K = k \right) f_K (k) 
%&=& \sum_{k = 0}^{\infty} P \left(\sum_{i = 0}^{k} Y^*_i \ge u \right) (1-\mu) \mu^k \\
= (1-\mu) \sum_{k = 1}^{\infty} P ( S_k^* \ge u) \mu^k.
\end{equation*} 
\qed

For example, suppose claims have a $\mbox{Geometric}(p)$ distribution
with mean $\mu=(1-p) /p$. The net profit condition $\mu<1$ implies $p>1/2$.
We have seen that the associated equilibrium distribution is again
$\mbox{Geometric}(p)$, and hence the $k$-th convolution is
$ \mbox{Negative Binomial}(k,p)$, $k\ge0$.
Straightforward calculations show that the
Pollaczeck--Khinchine formula gives the known solution for the probability of ruin,
\begin{equation}
\label{ruina-geo}
\psi (u)=\left(\frac{1-p} {p} \right)^{u + 1},\quad u\ge0.
\end{equation}

\begin{comment} *****
\begin{eqnarray}
\psi (u) &=& (1-\mu) \sum_{k = 1}^{\infty} P ( S_k^* \ge u) \, \mu^k \\
&=& (1-\mu) \sum_{k = 1}^{\infty} \sum_{i = 0}^{k-1} \binom {u + i-1} {i} \, p^i (1-p)^{u} \, \mu^k,\\
\nonumber &=& \left(1- \frac{1-p} {p} \right) \sum_{i = 0}^{\infty} \sum_{k = i + 1}^{\infty} \binom {u + i-1} {i} \, p^i (1-p)^{u} \, \left(\frac{1-p} {p} \right)^k \\
\nonumber &=& \left(1- \frac{1-p} {p} \right) \sum_{i = 0}^{\infty} \binom {u + i-1} {i} \, p^i (1-p)^{u} \, \sum_{k = i + 1}^{\infty} \left(\frac{1-p} {p} \right)^k \\
\nonumber &=& (1-p)^{u + 1} \sum_{i = 0}^{\infty} \binom {u + i-1} {i} \, (1-p)^i / p \\
\label{ruina-geo} &=& \left(\frac{1-p} {p} \right)^{u + 1}.
\end{eqnarray}
\end{comment} *****

This includes in the same formula the case $u=0$.
In the following section we will consider claims that have a mixture of some
distributions.

\begin{comment}
The following section describes and analyzes the negative binomial mixture distribution ( MIBN cambiar a NBMD ). One of the applications of equality (\ref{psi-geo-comp-d}) in this work is found in the proof of the Proposition \ref{psibn}, where the ruin probability is calculated when claims have MIBN distribution. Then, in Section \ref{sectionfive}, the Pollaczeck-Khinchine formula is applied within the proof of the main Theorem \ref{principal-myb-theorem}.
\end{comment}

%-----------------------------------------------------------------

\section{Negative binomial mixture}
\label{fivefour}

Negative binomial mixture (NBM) distributions 
will be used to approximate the ruin probability when claims
have a mixed Poisson (MP)
distribution. Although NBM distributions are the analogue of Erlang mixture distributions,
they cannot be used to approximate any discrete distribution with non-negative support. However, it turns out that they can approximate mixed Poisson
distributions. This is stated in \cite[Theorem 1] {1997steutel},
where the authors define NBM distributions those with probability
generating function
\begin{equation*}
G(z) = \lim_{m \rightarrow \infty} \sum_{k = 1}^m q_{k, m} \left(\frac{1-p_{k, m}} {1-p_{k , m} \, z} \right)^{r_{k, m}},\quad z<1,
\end{equation*}
where $ q_{k, m} $ are positive numbers and sum $1$ over index $k$.
This is a rather general definition for a NBM distribution.
In this work we will consider a particular case of it.\\

We will denote by $\mbox{nb}(k,p)(x)$ the probability function of a negative binomial distribution
with parameters $k$ and $p$, and by $\mbox{NB}(k,p)(x)$ its distribution function, namely,
for $x\ge0$,
\begin{equation*}
\mbox{nb}(k, p) (x) = \binom {k + x-1} {x} p^k (1-p)^x,\,\mbox{ and }\,\mbox{NB}(k, p) (x) = 1- \sum_{i = 0}^{k-1} \mbox{nb}(x + 1,1-p) (i).
\end{equation*}

\begin{definition}
Let $q_1, q_2,\ldots$ be a sequence of numbers such that $q_k\ge0$ and $\sum_{k = 1}^\infty q_k=1$.
A negative binomial mixture distribution with parameters $\boldsymbol{\pi}= (q_1, q_2, \ldots) $ and $p\in(0,1)$,
denoted by $\mbox{NBM}(\boldsymbol{\pi},p)$, is a discrete distribution with probability function
\begin{equation*}
f (x) = \sum_{k = 1}^\infty q_k \cdot \mbox{nb}(k, p) (x),\quad x \ge 0.
\end{equation*}
\end{definition}

It is useful to observe that any NBM distribution can be written as a compound sum
of geometric random variables. Indeed, let $N$ be a discrete random variable with
probability function $q_k=f_{N}(k)$, $k\ge1$, and define $S_N=\sum_{i = 1}^N X_i$,
where $ X_1, X_2, \ldots $ are i.i.d. r.v.s
$ \mbox{Geometric}(p) $ distributed and independent of $N$. Then
\begin{equation*}
\sum_{k = 1}^\infty q_k \cdot \mbox{nb}(k, p) (x)
=\sum_{k = 1}^\infty q_k\cdot P \left(\sum_{i = 1}^k X_i = x \right)
=P (S_N= x),\quad x\ge0.
\end{equation*}
Thus, given any $\mbox{NBM}(\boldsymbol{\pi},p)$ distribution with
$\boldsymbol{\pi} = (f_N (1), f_N (2), \ldots)$, we have the representation
\begin{equation}
\label{mibn-as-compound-sum}
S_N=\sum_{i = 1}^N X_i \sim \mbox{NBM}(\boldsymbol{\pi}, p).
\end{equation}
In particular,
\begin{eqnarray}
\label{532}
E ( S_N )&=&E (N) \left(\frac{1-p} {p} \right),\\
F_{S_N}(x) &=& \sum_{k = 1}^\infty f_{N}(k) \cdot \mbox{NB}(k, p) (x),\quad x\ge0,
\end{eqnarray}
and the p.g.f. has the form $G_{S_N}(r)=G_N( G_X(r))$.
The following is a particular way to write the distribution function of a NBM distribution.

\begin{proposition}
\label{F-de-mibn}
Let $ S_N \sim \mbox{NBM}(\boldsymbol{\pi}, p) $, where
$\boldsymbol{\pi} = (f_N (1), f_N (2 ), \ldots) $ for some discrete r.v. $N$.
For each $x\ge0$, let $Z\sim\mbox{NegBin}(x + 1,1-p)$. Then 
\begin{equation}\label{533}
F_{S_N}(x) = E (F_N (Z)),\quad x\ge0.
\end{equation}
\end{proposition}
\textbf{Proof.}
\begin{eqnarray*}
F_{S_N}(x) 
&=& \sum_{k = 1}^\infty f_{N}(k) \cdot \mbox{NB}(k, p) (x) \\
           &=& \sum_{k = 1}^\infty f_{N}(k) \, \left [1- \sum_{i = 0}^{k-1} \mbox{nb}(x + 1, 1-p) (i) \right] \\
%           &=& \sum_{k = 1}^\infty f_{N}(k) \, \sum_{i = k}^{\infty} \mbox{nb}(x + 1,1-p) (i ) \\
           &=& \sum_{i = 0}^{\infty} \left[\sum_{k = 1}^i f_{N}(k)\right] \, \mbox{nb}(x + 1,1-p) (i) \\
                      &=& E ( F_{N}(Z)).
\end{eqnarray*}
\qed

We will show next that the equilibrium distribution associated to a NBM distribution
is again NBM. For a distribution function $F(x)$, $\overline{F}(x)$ denotes $1-F(x)$.

\begin{proposition}\label{mibn-balance}
Let $S_N \sim \mbox{NBM}(\boldsymbol{\pi}, p) $, with
$\boldsymbol{\pi} = (f_N (1), f_N (2 ), \ldots) $ and $E(N)<\infty$.
The equilibrium distribution of $ S_N $ is $ \mbox{NBM}(\boldsymbol{\pi}_e, p) $, where
$\boldsymbol{\pi}_e = (f_{Ne}(1), f_{Ne}(2), \ldots) $ and
\begin{equation}
\label{534}
f_{Ne}(j) = \overline{F}_{N}(j-1) / E (N),\quad j \ge 0.
\end{equation}
\end{proposition}
\textbf{Proof.}
%The definition of $ f_e $ is used and then the Proposition \ref{F-de-mibn} is applied .
\begin{equation*}
f_e (x) = \frac{\overline{F}_{S_N}(x)}{E(S_N)}=\frac{p\sum_{i = 0}^{\infty} \overline{F}_N (i) \binom {x + i} {i} p^i (1-p)^{x + 1}} {(1-p) E (N)}=\sum_{i = 0}^{\infty} \frac{\overline{F}_N (i)} {E (N)} \binom {x + i} {i} p^{i + 1}(1-p)^x.
\end{equation*}
%\begin{eqnarray*}
%f_e (x) &=& \frac{1} {E ( S_N )}(1- F_{S_N}(x)) \\
%      &=& \frac{p} {(1-p) E (N)} \sum_{i = 0}^{\infty} \overline{F}_N (i) \binom {x + i} {i} p^i (1-p)^{x + 1} \\
%      &=& \sum_{i = 0}^{\infty} \frac{\overline{F}_N (i)} %{E (N)} \binom {x + i} {i} p^{i + 1}(1-p) x,
%\end{eqnarray*}
Naming $ j = i + 1 $,
\begin{equation*}
f_e (x) = \sum_{j = 1}^{\infty} \frac{\overline{F}_{N}(j-1)} {E (N)} \binom {j + x- 1} {x} p^{j}(1-p)^{x} 
     = \sum_{j = 1}^{\infty} f_{Ne}(j) \cdot \mbox{nb}(j, p) (x).
\end{equation*}
\qed

It can be checked that~(\ref{534}) is a probability function. It is the equilibrium
distribution associated to $N$. In what follows,
a truncated geometric distribution will be used. This is denoted by
$\mbox{TGeometric}(\rho)$, where $0<\rho<1$, and defined by the probability function
$ f(k) = \rho (1- \rho)^{k-1} $, for $ k \ge 1 $.\\

The following proposition states that a compound geometric  NBM distribution is again NBM.
This result is essential to calculate the ruin probability when claims have NBM distribution.

\begin{proposition}
\label{sum-compound-mibn}
Let $ M \sim \mbox{TGeometric}(\rho)$ and let $N_1,N_2,\ldots$ be a sequence of independent
random variables with identical distribution $ \boldsymbol{\pi} = (f_N (1), f_N (2), \ldots) $.
Let $ S_{N_1}, S_{N_2}, \ldots $ be random variables with $ \mbox{NBM}(\boldsymbol{\pi}, p)$
distribution. Then
\begin{equation}
S:=\sum_{j = 1}^{M} S_{N_j} \sim \mbox{NBM}(\boldsymbol{\pi}^*, p),
\end{equation}
where $\boldsymbol{\pi}^*=(f_{N^*} (1), f_{N^*} (2), \ldots)$ is the distribution of
$ N^* = \sum_{j = 1}^{M} N_j $ and is  given by
\begin{eqnarray}
\label{536}
f_{N^*}(1) &=& \rho \, f_N (1), \\
\label{537}
f_{N^*}(k) &=& (1- \rho) \sum_{i = 1}^{k-1} f_N (i) \, f_{N^*}(k-i) + \rho \, f_{N}(k),\quad k\ge 2.
\end{eqnarray}
\end{proposition}

\textbf{Proof.} For $ x \ge 1 $ and $ m \ge 1 $,
\begin{equation}\label{aux1} 
P (S = x \mid M = m) =  P  \left(\sum_{j = 1}^{m} S_{N_j} = x \right)= P \left(\sum_{j = 1}^{m} \sum_{i = 1}^{N_j} X_{i \, j} = x \right) =P \left(\sum_{\ell = 1}^{N_m} X_{\ell} = x \right),             
\end{equation}
where $ N_m = \sum_{i = 1}^{m} N_i $ and $ X_{\ell} \sim \mbox{Geometric}(p) $ for $ \ell \ge 1 $. Therefore,
\begin{equation*}
P (S = x) = \sum_{m = 1}^{\infty} P (S = x \mid M = m) f_M (m) = \sum_{m = 1}^{\infty} P \left(\sum_{\ell = 1}^{N_m} X_{l} = x \right) f_M (m) = P \left(\sum_{\ell = 1}^{N^*} X_{\ell} = x \right),
\end{equation*}
where $ N^* = \sum_{j = 1}^{M} N_j $. Using Panjer's formula it can be shown that $ N^*$
has distribution $ \boldsymbol{\pi}^* $ given by~(\ref{536}) and~(\ref{537}).
Since $ X_{\ell} \sim \mbox{Geometric}(p) $,
$ \sum_{\ell = 1}^{N^*} X_{\ell} \sim \mbox{NBM}(\boldsymbol{\pi}^*, p) $.
\begin{comment} ***
\textcolor{red}{----------------------------------------------}

\textcolor{red}{For any random variable $X$ we will denote its probability generating function by $G_X(r)=\sum_{i=0}^\infty r^i P(X=r)$. Let $ S = \sum_{j = 1}^{M} S_{N_j} $, $ M \sim \mbox{TGeometric}(\rho) $ and $ S_{N_j} = \sum_{i = 1}^{N_j} X_{i \, j} $, where $ P (N_j = k) = f_{N}(k) $ and $ X_{i \, j} \sim \mbox{Geometric}(p) $. }
\end{comment}
Lastly, we consider the probability of the event $(S=0)$.
\begin{equation*}
P (S = 0) = \sum_{k = 1}^{\infty} f_{N^*}(k) \, \mbox{nb}(k, p) (0) = \sum_{k = 1}^{\infty} f_{N^*}(k) \, p^k 
= f_{N^*}(1) \, p + \sum_{k = 2}^{\infty} f_{N^*}(k) \, p^k.
\end{equation*}
Substituting $ f_{N^*}(k) $ from (\ref{536}) and~(\ref{537}), one obtains
\begin{eqnarray*}
P (S = 0) 
%&=& \rho \, f_N (1) \, p + \sum_{k = 2}^{\infty} \left [(1- \rho) \sum_{i = 1}^{k-1} f_N (i) \, f_{N^*}(k-i) + \rho \, f_{N}(k) \right] \, p^k \\
%      &=& \rho \, f_N (1) \, p + (1- \rho) \, \sum_{k = 2}^{\infty} \sum_{i = 1}^{k-1} f_N (i ) \, f_{N^*}(k-i) \, p^k + \rho \, \sum_{k = 2}^{\infty} f_{N}(k) \, p^k \\
%      &=& \rho \, \sum_{k = 1}^{\infty} f_{N}(k) \, p^k + (1- \rho) \, \sum_{i = 1}^{\infty} f_N (i) \, \sum_{k = i + 1}^{\infty} f_{N^*}(k-i) \, p^k,\\
% &=& \rho \, G_N (p) + (1- \rho) \, \sum_{i = 1}^{\infty} f_N (i) \, \sum_{j = 1}^{\infty} f_{N^*}(j) \, p^{i + j} \\
%      &=& \rho \, G_N (p) + (1- \rho) \, \sum_{i = 1}^{\infty} f_N (i) \, p^{i} \, \sum_{j = 1}^{\infty} f_{N^*}(j) \, p^{j} \\
      &=& \rho \, G_N (p) + (1- \rho) \, G_N (p) \, P (S = 0).     
\end{eqnarray*}
Therefore,
\begin{equation}
\label{PSsea0}
P (S = 0)=\frac{\rho \, G_N (p)} {1- (1- \rho) \, G_N (p)}= G_M ( G_N (p))= G_M ( G_N ( G_{X_{i \, j}}(0))).
\end{equation}
%On the other hand, assuming that $ S \sim \mbox{NBM}(\boldsymbol{\pi}^*, p) $,
%\begin{equation}\label{PSsea0}
%P (S = 0) = G_S (0) = G_M ( G_N ( G_{X_{i \, j}}(0))) = G_M ( G_N (p)) = \frac{\rho \, G_N (p)} {1- (1- \rho) \, G_N (p)}.
%\end{equation}
The last term is the p.g.f. of a $\mbox{NBM}(\boldsymbol{\pi}^*, p)$ distribution evaluated
at zero.
\qed

From~(\ref{536}) and (\ref{537}), it is not difficult to derive a recursive formula for
$\overline{F}_{N^*}(k)$, namely,

\begin{comment}
a formula can be derived to calculate $ \overline{F}_{N^*}(k) = \sum_{i = k + 1}^\infty f_{N^*}(i) $. First, for $ k = 1 $,
\begin{equation*}
F_{N^*}(1) = f_{N^*}(1) = \rho \, f_N (1),
\end{equation*}
therefore, $ \overline{F}_{N^*}(1) = 1- \rho \, f_N (1) $. For $ k \ge 2 $,
\begin{eqnarray*}
F_{N^*}(k) &=& \sum_{i = 1}^{k} f_{N^*}(i) = f_{N^*}(1) + \sum_{i = 2}^{k} \left [(1- \rho) \sum_{j = 1}^{i-1} f_N (j) \, f_{N^*}(i-j) + \rho \, f_{N}(i) \right] \\
   &=& \rho \, f_N (1) + (1- \rho) \sum_{i = 2}^{k} \sum_{j = 1}^{i-1} f_N (j) \, f_{N^*}(i-j) + \rho \, \sum_{i = 2}^{k} f_{N}(i) \\
   &=& \rho \, \sum_{i = 1}^{k} \, f_{N}(i) + (1- \rho) \sum_{j = 1}^{k-1} f_N (j ) \, \sum_{i = j + 1}^{k} f_{N^*}(i-j),
\end{eqnarray*}
making the change of index $ m = i-j $, the following expression is reached,
\begin{equation*}
F_{N^*}(k) =\rho \, F_N (k) + (1- \rho) \sum_{j = 1}^{k-1} f_N (j) \, F_{N^*}(k-j).    
\end{equation*}
Then,
\end{comment}
\begin{equation}
\label{539}
\overline{F}_{N^*}(k) = (1 - \rho) \, \sum_{j = 1}^k f_N (j) \, \overline{F}_{N^*}( k-j) + \overline{F}_N (k),\quad k\ge 1.
\end{equation}

The following result establishes a formula to calculate the ruin probability when claims have a
NBM distribution.

\begin{theorem}
\label{psibn}
Consider the Gerber-Dickson model with claims having a
$\mbox{NBM}(\boldsymbol{\pi}, p)$ distribution, where $\boldsymbol{\pi}=(f_N (1), f_N (2),\ldots)$
and $E (N)<\infty$. For $u\ge1$ define $Z_u\sim\mbox{NegBin}(u, 1-p)$.
Then the ruin probability can be written as
\begin{equation}
\label{psiu-mibn}
\psi (u) = \sum_{k = 0}^{\infty} \overline{C}_k \cdot P ( Z_u = k) = E (\overline{C}_{Z_u}),\quad u\ge 1,
\end{equation}
%\begin{eqnarray} \label{psiu-mibn}
%\psi (u) &=& \sum_{k = 0}^{\infty} \overline{C}_k \, P ( %Z_u = k) \\
%\label{psiu-mibn-as-hope}       
%&=& E (\overline{C}_{Z_u}),\quad u \ge 1 ,
%\end{eqnarray}
where the sequence $\left\{\overline{C}_k \right\}_{k = 0}^\infty $ is given by
\begin{eqnarray}
\label{543} \overline{C}_0 &=& E (N) (1-p) / p, \\
\label{544} \overline{C}_k &=& \overline{C}_0 \, \left [\sum_{i = 1}^{k} f_{Ne}(i) \, \overline{C}_{k-i} + \overline{F}_{Ne}(k) \right],\quad k\ge1, \\
\label{545} f_{Ne}(i) &=& \frac{\overline{F}_N (i-1)} {E (N)},\quad i\ge1.
\end{eqnarray}
\end{theorem}
\textbf{Proof.}
Let $R_0 = \sum_{j = 1}^{M_0} Y_{e, j}$,
where $ M_0 \sim \mbox{Geometric}(\rho) $ with $\rho=1-\psi(0)$,
and let $ Y_{e, 1}, Y_{e, 2}, \ldots $ be r.v.s distributed according to the equilibrium distribution
associated to $ \mbox{NBM}(\boldsymbol{\pi}, p) $ claims. By Proposition \ref{mibn-balance},
we know this equilibrium distribution is $\mbox{NBM}(\boldsymbol{\pi}_e, p)$,
where $\boldsymbol{\pi}_e$ is given by $ f_{Ne}(j) = \overline{F}_N (j-1 ) / E (N) $, $j\ge1$.
By~(\ref{psi-geo-comp-d}), for $ u \ge 1 $,
\begin{eqnarray*}
\psi (u) &=& P ( R_0 \ge u) \\
       &=& P ( R_0 \ge u \mid M_0 > 0) P ( M_0 > 0) + P ( R_0 \ge u \mid M_0 = 0) P ( M_0 = 0) \\
%       &=& P ( R_0 \ge u \mid M_0 > 0) \, (1- \rho) \\
       &=& (1- \rho) \, P (R \ge u),
\end{eqnarray*}
where $R\sim\sum_{j = 1}^{M} Y_{e, j}$ with
$ M \sim \mbox{TGeometric}(\rho) $ with probability function $ f_M (k) = \rho (1- \rho)^{k-1} $, for $ k \ge 1 $.
%and by the Proposition \ref{num-de-jumps-d}, $ \rho = 1 - \psi (0) $.
By Proposition \ref{sum-compound-mibn}, $ R \sim \mbox{NBM}(\boldsymbol{\pi}^*, p) $, where
$\boldsymbol{\pi}^* $ is given by equations (\ref{536}) and (\ref{537}).
Now define
\begin{equation}\label{C-como-proba-de-cola}
\overline{C}_k = (1- \rho) \overline{F}_{N^*}(k),\quad k\ge0.
\end{equation}
Therefore, using (\ref{533}),
\begin{equation*}
\psi (u)=
(1- \rho) P (R> u) = (1- \rho) E \left(\overline{F}_{N^*}( Z_{u}) \right) = \sum_{k = 0}^\infty \overline{C}_k P ( Z_{u} = k).
\end{equation*}
Finally, we calculate the coefficients $\overline{C}_k $ where $ \rho = 1- \psi (0) = 1-E (N) (1-p) / p $.
First,
\begin{equation*}
\overline{C}_0 = (1- \rho) \overline{F}_{N^*}(0) = 1- \rho = E (N) (1-p) / p,
\end{equation*}
and by (\ref{539}),
\begin{equation*}
\overline{C}_k = (1- \rho) \overline{F}_{N^*}(k) = \overline{C}_0 \, \left [\, \sum_{i = 1}^k f_{Ne}(i) \overline{C}_{k-i} + \overline{F}_{Ne}(k) \, \right],\quad  k\ge1.
\end{equation*}
\qed

As an example consider claims with a geometric distribution. This is a NBM distribution with
$\boldsymbol{\pi}=(1,0,0, \ldots)$. Equations~(\ref{543}--\ref{545}) yield
\begin{displaymath}
\overline{C}_k=\left((1-p) / p \right)^{k + 1},\quad k\ge0.
\end{displaymath}
\begin{comment}
\begin{eqnarray*}
\overline{C}_0 &=& (1-p) / p. \\
\overline{C}_1 &=& (1-p) / p \left( f_{Ne}(1) (1-p) / p + \overline{F}_{Ne}(1) \right) = \left((1-p) / p \right)^2. \\
\overline{C}_2 &=& (1-p) / p \left( f_{Ne}(1) ((1-p) / p)^2 + f_{Ne}(2) (1-p) / p + \overline{F}_{Ne}(2) \right) = \left((1-p) / p \right)^3. \\
& \vdots & \\
\overline{C}_k &=& \left((1-p) / p \right)^{k + 1}.
\end{eqnarray*}
\end{comment}
Substituting in~(\ref{psiu-mibn}) together with $\psi(0)=(1-p)/p$, we recover the known
solution $\psi (u)=\left((1-p)/p\right)^{u + 1},\quad u\ge0$.

%%%%%%%%%%%%%%%%%%%%%%%%%%%%%%%%%%%%%%%%%%%%%%%%%%%%%
\section{Mixed Poisson}
\label{cincocinco}

This section contains the definition of a mixed Poisson distribution and some of its
relations with NBM distributions.

\begin{definition}
\label{definition-Poisson-mixed}
Let $X$ and $\Lambda$ two non-negative random variables. If
$ X \mid (\Lambda = \lambda) \sim \mbox{Poisson}(\lambda) $, then we say that $X$ has
a mixed Poisson distribution with mixing distribution $ F_{\Lambda}$.
In this case, we write $ X \sim \mbox{MP}(F_\Lambda) $.
\end{definition}

Observe the distribution of $ X \mid (\Lambda = \lambda)$
is required to be Poisson, but
the unconditional distribution of $X$, although discrete, is not
necessarily Poisson. A large number of examples of these distributions can be found
in~\cite{2005karlis} and a study of their general properties is given in~\cite{1997grandell}.
In particular, it is not difficult to see that $E(X)=E(\Lambda)$ and the p.g.f. of $X$ can be
written as
\begin{equation}
\label{fgp-poisson-mixed}
G_X (r) = \int_0^\infty e^{- \lambda (1-r)} dF_{\Lambda}(\lambda),\quad r<1.
\end{equation}

The following proposition establishes a relationship between the Erlang mixture distribution and the negative binomial distribution. The former will be denoted by
$\mbox{ErlangM}(\boldsymbol{\pi},\beta)$, with similar meaning for the parameters as in the
notation $\mbox{NBM}(\boldsymbol{\pi},p)$ used before.
In the ensuing calculations the probability function of a $\mbox{Poisson}(\lambda)$ distribution
will be denoted by $\mbox{poisson}(\lambda)(x)$.

\begin{proposition}
\label{MIENBM}
Let $\Lambda$ be a random variable with distribution
$\mbox{ErlangM}(\boldsymbol{\pi}, \beta) $.
The distributions $\mbox{MP}(F_\Lambda)$ and $ \mbox{NBM}(\boldsymbol{\pi}, \beta / (\beta +1))$ are the same.
\end{proposition}
\textbf{Proof.} 
Let $ X \sim \mbox{MP}(F_\Lambda) $. For $x\ge0$,
\begin{eqnarray*}
P (X = x) &=& \int_0^\infty poisson(\lambda)(x) \cdot \sum_{k = 1}^\infty q_k \cdot
\mbox{erl}(k, \beta) (\lambda)\, d\lambda \\
&=& \sum_{k = 1}^\infty q_k \cdot
\left(\frac{\beta} {\beta + 1} \right)^k \left(\frac{1} {\beta + 1} \right)^x \frac{(k + x-1)!} {(k-1)! \, x!} \\
 &=& \sum_{k = 1}^\infty q_k \cdot \mbox{nb}(k, \beta / (\beta + 1)) (x).
\end{eqnarray*}
\qed

As a simple example consider the case $\Lambda \sim \mbox{Exp}(\beta)$ and
$\boldsymbol{\pi}=(1,0,0,\ldots)$. By Proposition \ref{MIENBM},
$P (X = x) = \mbox{nb}(1, \beta / (\beta + 1)) (x)$ for $x\ge0$.
That is, $ X \sim \mbox{Geometric}(p) $ with $ p = \beta / (\beta + 1) $.\\

Next proposition will be useful to show that
a MP distribution can be approximated by NBM distributions.  Its proof can be
found in~\cite{1997grandell}.

\begin{proposition}\label{proposition-grandell}
Let $ \Lambda_1, \Lambda_2, \ldots $ be positive random variables with distribution functions
$ F_1, F_2, \ldots $ and let $ X_1, X_2, \ldots $ be random variables such that
$ X_i \sim \mbox{MP}(F_i) $, $ i \ge 1 $. Then $ X_n \xrightarrow {D} ​​X $, if and only if, $ \Lambda_n \xrightarrow {D} ​​\Lambda $, where $ X \sim \mbox{MP}(F_\Lambda) $.
\end{proposition}
%As usual, the notation $ Z_n \xrightarrow {D} ​Z $ means that the sequence $ Z_1, Z_2, \ldots $ converges in distribution to $ Z $. 

Finally we establish how to approximate an MP distribution.

\begin{proposition}\label{propos-PXn-a-PX}
Let $ X \sim \mbox{MP}(F_\Lambda) $, and let $ X_n $ be a random variable with distribution $ \mbox{NBM}(\boldsymbol{\pi}_n, p_n) $ for $ n \ge 1 $, where $p_n = n / (n + 1)$, $\boldsymbol{\pi}_n = (q (1, n), q (2, n), \ldots)$ and $q (k, n) = F_{\Lambda}(k / n) - F_{\Lambda}((k-1) / n)$.
Then $ X_n \xrightarrow {D} ​​X $.
\end{proposition}

\textbf{Proof.}
First, suppose that $ F_\Lambda $ is continuous. Let $ \Lambda_1, \Lambda_2, \ldots $ be random variables, where $\Lambda_n$ has distribution given by the following Erlangs mixture (see \cite{2016santana}),
\begin{equation}
\label{efene-dd}
F_{n}(x) = \sum_{k = 1}^{\infty} q (k, n) \cdot \mbox{Erl}(k, n) (x),\quad x>0,
\end{equation}
with $ q (k, n) = F_\Lambda (k / n) -F_\Lambda ((k-1) / n) $.
It is known~\cite{2016santana} that
\begin{equation*}
\lim \limits_{n \rightarrow \infty} F_{n}(x) = F_{\Lambda}(x), \quad x>0.
\end{equation*}
Then, by Proposition~\ref{proposition-grandell}, $ X_n \xrightarrow {D} ​​X $,
where $X_n\sim\mbox{MP}(F_n)$. This is an
$\mbox{NBM}(\boldsymbol{\pi}, p_n) $ by Proposition~\ref{MIENBM} where
$ \boldsymbol{\pi} = (q (1, n), q (2, n), \ldots) $ and $ p_n = n / (n + 1) $.\\

Now suppose $ F_\Lambda $ is discrete.
Let $ Y_n \sim \mbox{NegBin}(\lambda n, n / (n + 1)) $, where $\lambda$ and $ n $ are positive integers and let $ Z \sim\mbox{Poisson}(\lambda) $. The probability generating functions of these random
variables satisfy
\begin{equation*}
\lim_{n \rightarrow \infty} G_{Y_n}(r) = \lim_{n \rightarrow \infty} \left(1+ \frac{1-r} {n} \right)^{- \lambda n}=\exp\{-\lambda (1-r)\} = G_Z (r).
\end{equation*}
Thus,
\begin{equation}\label{chida}
Y_n \xrightarrow {D} ​​Z.
\end{equation}
On the other hand, suppose that $ X $ is a mixed Poisson random variable with probability function $ f_X (x) $ for $ x \ge 0$ and mixing distribution $ F_\Lambda (\lambda) $ for $ \lambda \ge 1$. Let $ \{X_n \}_{n = 1}^\infty $ be a sequence of random variables with distribution
\begin{equation}\label{efene-d}
f_n (x) = \sum_{k = 1}^{\infty} q (k, n) \cdot
\mbox{nb} \left(k, \frac{n}{n+1}\right) (x),\quad n\ge 1,\, x\ge 0,
\end{equation}
where $ q (k, n) = F_\Lambda (k / n) -F_\Lambda ((k-1) / n) $. Note that for any natural number $ n $, if $ k $ is not a multiple of $ n $, then $ q (k, n) = 0 $. Let $ k = \lambda \, n $ with $ \lambda \ge 1$.
Then $q (k, n) = F_\Lambda (\lambda)-F_\Lambda (\lambda-1/n)=f_\Lambda (\lambda)$. Therefore, for $ x\ge 0$,
\begin{displaymath}
f_{X_n}(x)
=\sum_{\lambda = 1}^\infty q (\lambda \, n, n) \cdot \mbox{nb}(\lambda \, n, n / (n + 1)) (x)
=\sum_{\lambda = 1}^\infty f_\Lambda (\lambda) \cdot \mbox{nb}(\lambda \, n, n / (n + 1)) (x).
\end{displaymath}
Therefore,
\begin{displaymath}
\lim_{n \rightarrow \infty} f_{X_n}(x) 
= \sum_{\lambda = 1}^\infty
f_\Lambda (\lambda) \cdot \lim_{n \rightarrow \infty} \mbox{nb}(\lambda \, n, n / (n + 1)) (x)
= \sum_{\lambda = 1}^\infty
f_\Lambda (\lambda) \cdot \mbox{poisson}(\lambda) (x).
\end{displaymath}
\qed

For $ X_n\sim\mbox{NBM}(\boldsymbol{\pi}_n, p_n) $ as in the previous statement, it easy to see that
\begin{equation}
\label{media-exn-finite}
E (X_n) <1.
\end{equation}
As a consequence of Proposition \ref{propos-PXn-a-PX}, for $ X \sim \mbox{MP}(F_{\Lambda}) $, its probability function can be approximated by NBM distributions with
suitable parameters. That is, for sufficiently large values of $ n $,
\begin{equation}\label{NBM-approximants-aX}
P (X = x) \approx \sum_{k = 1}^{\infty} q (k, n) \cdot \mbox{nb}(k, p_n) (x),
\end{equation}
where $q (k, n) = F_{\Lambda} \left(k / n \right) - F_{\Lambda} \left((k-1) / n \right)$ and
$p_n=n/(n + 1)$.

%-----------------------------------------------------------------

\section{Ruin probability approximations}
\label{sectionfive}

We here consider the case when claims in the Gerber-Dickson risk model have distribution
function $F\sim\mbox{MP}(F_{\Lambda})$.
Let $\psi_n(u)$ denote the ruin probability when claims have distribution $F_n(x)$ as defined
in Proposition~\ref{propos-PXn-a-PX}.
If $n$ is large enough, $ F_n (x) $ is close to $ F (x) $, and is expected that
$ \psi_n (u) $ will be close to $ \psi (u ) $, the unknown ruin probability.
This procedure is formalized in the following theorem.

\begin{theorem} \label{principal-myb-theorem}
If claims in the Gerber-Dickson model have a $\mbox{MP}(F_{\Lambda})$ distribution, then
\begin{equation*}
\psi(u) = \lim \limits_{n \rightarrow \infty} \psi_n (u),\quad u\ge0,
\end{equation*}
where
\begin{equation}\label{560}
\psi_n (u) =\sum_{k = 0}^{\infty} \overline{C}_{k, n} \, P (Z = k) = E \left(\overline{C}_{Z, n} \right),
\end{equation}
%\begin{eqnarray}
%\label{560} \psi_n (u) &=& \sum_{k = 0}^{\infty} \overline{C}_{k, n} \, P (Z = k) \\
%\label{561} &=& E \left(\overline{C}_{Z, n} %\right).
%\end{eqnarray}
with $ Z\sim \mbox{NegBin}(u,1/(1 + n))$.
The sequence $ \left\{\overline{C}_{k, n} \right\}_{k = 0}^\infty $ is determined by
\begin{eqnarray}
\label{562} \overline{C}_{0, n} &=& \sum_{j = 0}^\infty \overline{F}_\Lambda (j / n) / n, \\
\label{563} \overline{C}_{k, n} &=& \overline{C}_{0, n} \left [\sum_{i = 1}^{k} f_{Ne}(i) \overline{C}_{k-i, n} + \overline{F}_{Ne}(k) \right],\quad k\ge1, \\
\label{564} f_{Ne}(i) &=& \frac{\overline{F}_\Lambda ((i-1) / n)} {\sum_{j = 0}^\infty \overline{F}_\Lambda (j / n)},\quad i\ge1.
\end{eqnarray}
\end{theorem}

\textbf{Proof.} Suppose $ X \sim \mbox{MP}(F_\Lambda) $ with $ E (X) <1 $ and
equilibrium probability function $f_e(x)$.
%and let $f_e(x)$ and $G_e(r)$ be the related equilibrium probability function and generating probability function, respectively.  
Let $ X_1, X_2, \ldots $ be an approximating sequence of $ \mbox{NBM}(\boldsymbol{\pi}_n , p_n) $
random variables to $X$, where $\boldsymbol{\pi}_n = (q (1, n), q (2, n), \ldots)$, with
$q (k, n) = F_{\Lambda}(k / n) - F_{\Lambda}((k-1) / n)$ and $p_n = n / (n + 1)$. That is,
\begin{equation}
f_{X_n}(x) = \sum_{k = 1}^\infty q (k, n) \cdot \mbox{nb}(k, n / (n + 1)) (x),\quad x\ge 0.
\end{equation}
%and probability generating función, respectively. 
By (\ref{532}),
\begin{equation*}
E ( X_n ) =  \sum_{k = 1}^\infty k \, q (k, n)\cdot\frac{1/(n+1)}{n/(n+1)}
= \sum_{k = 1}^\infty (k/n) \cdot [\,F_{\Lambda}(k / n) - F_{\Lambda}((k-1) / n)\,].
\end{equation*}
Taking the limit,
\begin{equation}\label{aaa3}
\lim_{n \rightarrow \infty} E ( X_n )=\int_0^\infty x\,dF_{\Lambda}(x)= E (\Lambda)= E (X).
\end{equation}
Now, by Proposition \ref{propos-PXn-a-PX}, since $ X_n \xrightarrow {D} ​​X $, we have
\begin{equation*}
\lim_{n \rightarrow \infty} \overline{F}_{X_n}(x) = \overline{F}_X (x),\quad x\ge 0.
\end{equation*}
Combining the above with (\ref{aaa3}),
\begin{equation*}
\lim_{n \rightarrow \infty} \frac{\overline{F}_{X_n}(x)} {E (X_n)} = \frac{\overline{F}_X (x)} {E (X)}.
\end{equation*}
This means the equilibrium probability function $f_{e, n}(x)$ associated to $f_{X_n}(x)$
satisfies
\begin{equation}\label{aaa4}
\lim_{n \rightarrow \infty} f_{e, n}(x) = f_e (x),\quad x \ge 0.
\end{equation}
%%%%%%%%%%%%%%%%%%%%%%%%%%%%%%%%%%%%%%%%%
\begin{comment}
The probability generating function of $S_{k, n}: = \sum_{i = 1}^k X_{ni}$ is,
\begin{equation*}
G_{S_{k, n}}(r) = (G_n (r))^k,\quad r <1 .
\end{equation*}
Using (\ref{aaa4}) and due to the continuity of the probability generating function,
\begin{equation}\label{aaa5}
\lim_{n \rightarrow \infty} G_{S_{k, n}}(r) = \lim_{n \rightarrow \infty}( G_n (r))^k = (\lim_{n \rightarrow \infty} G_n (r))^k = ( G_e (r))^k,
\end{equation}
On the other hand, let  $S_k: = \sum_{i = 1}^k X_{e, i}$, where $ X_{e, 1}, X_{e, 2}, \ldots $ are i.i.d.r.v. with probability function $ f_e (x) $, then his probability generating function is,
\begin{equation}\label{568}
G_{S_k}(r) = (G_e (r))^k.
\end{equation}
Equaling (\ref{568}) and (\ref{aaa5}),
\begin{equation*}
\lim_{n \rightarrow \infty} G_{S_{k, n}}(r) = G_{S_k}(r),\quad r <1.
\end{equation*}
Thus,
\end{comment}
%%%%%%%%%%%%%%%%%%%%%%%%%%%%%%
Using probability generating functions and (\ref{aaa4}), it is also easy to show that
for any $k\ge1$,
\begin{equation}\label{aaa6}
\lim_{n \rightarrow \infty} F_{e, n}^{* k}(x) = F_e^{* k}(x),\quad x\ge0.
\end{equation}
Now, let $ X_{n1}, X_{n2}, \ldots $ be i.i.d. random variables with
probability function $ f_{e, n}(x) $ and set $S_{k, n}: = \sum_{i = 1}^k X_{ni}$.
By the Pollaczeck-Khinchine formula, for $u\ge0$,
\begin{equation*}
\psi_n (u) = \sum_{k = 1}^{\infty} P ( S_{k, n} \ge u) \, (1-E ( X_n )) E^k ( X_n ) = \sum_{k = 1}^{\infty}(1- F_{e, n}^{* k}(u-1)) \, (1-E ( X_n )) E^k( X_n ).
\end{equation*}
%\begin{eqnarray*}
%\psi_n (u) &=& \sum_{k = 1}^{\infty} P ( S_{k, n} \ge u) %\, (1-E ( X_n )) E ( X_n )^k \\
%         &=& \sum_{k = 1}^{\infty}(1-P ( S_{k, n} \le %u-1)) \, (1-E ( X_n )) E ( X_n )^k \\
%         &=& \sum_{k = 1}^{\infty}(1- F_{e, n}^{* k}(u-1)) %\, (1-E ( X_n )) E ( X_n )^k
%\end{eqnarray*}
Taking the $n\to\infty$ limit and using (\ref{aaa3}) and (\ref{aaa6}),
\begin{displaymath}
\lim_{n \rightarrow \infty} \psi_n (u)
%&=& \sum_{k = 1}^{\infty} \lim_{n \rightarrow \infty}(1- F_{e, n}^{* k}(u-1)) \, (1-E ( X_n )) E ( X_n )^k \\
= \sum_{k = 1}^{\infty}(1- F_{e}^{* k}(u-1)) \, (1-E (X)) E^k (X)
= \psi (u),\quad u\ge1.
\end{displaymath}
%The limit enters into the sum because of the Monotone Convergence Theorem for series. By (\ref{media-exn-finite}) it is true that $ E (X_n) <1 $, so $ E (X_n)^k $ is used as the dominant sequence, and $ \sum_{k = 1}^\infty E (X_n)^k = E (X_n) / (1-E (X_n)) <\infty $.\\ 
On the other hand, since claims $ X_n $ have a $\mbox{NBM}(\boldsymbol{\pi}_n, p_n) $ distribution, with
$\boldsymbol{\pi}_n = (q (1, n), q (2, n), \ldots)$, $\quad q (k, n) = F_{\Lambda}(k / n) - F_{\Lambda}((k-1) / n)$ and $p_n    = n / (n + 1)$, by Theorem~\ref{psibn},
\begin{equation*}
\psi_n (u) = \sum_{k = 0}^{\infty} \overline{C}_{k, n} \cdot P (Z = k) = E (\overline{C}_{Z, n}),\quad u\ge 1,
\end{equation*}
where $ Z \sim \mbox{NegBin}(u, 1 / (n + 1)) $ and the sequence $ \left\{\overline{C}_{k, n} \right\}_{k = 0}^\infty $ is given by

\begin{eqnarray*}
\overline{C}_{0, n} &=& E ( N_n ) / n, \\
\overline{C}_{k, n} &=& \overline{C}_{0, n} \, \left [\sum_{i = 1}^{k} f_{Ne}(i) \, \overline{C}_{k-i, n} + \overline{F}_{Ne}(k) \right],\quad k\ge1, \\
f_{Ne}(i) &=& \frac{\overline{F}_{N_n}(i-1)} {E ( N_n )},\quad i\ge1,
\end{eqnarray*}
where $N_n$ is the r.v. related to probabilities $q(k,n)$. Thus, it only remains to calculate the form of $ E (N_n) $ and $ \overline{F}_{N_n}(i-1) $.
%\begin{eqnarray*}
%E ( N_n ) &=& \sum_{j = 1}^\infty P ( N_n > j-1) \\
% &=& \sum_{j = 1}^\infty \sum_{i = j}^\infty P (N_n = i) \\
%       &=& \sum_{j = 1}^\infty \sum_{i = j}^\infty q (i, n) \\
 %       &=& \sum_{j = 1}^\infty \sum_{i = j}^\infty ( F_\Lambda (i / n) - F_\Lambda ((i-1) / n)) \\
% &=& \sum_{j = 1}^\infty (1-F_\Lambda ((j-1) / n)) \\
%        &=& \sum_{j = 0}^\infty \overline{F}_\Lambda (j / n).
%\end{eqnarray*}
\begin{equation*}
E ( N_n ) = \sum_{j = 1}^\infty P ( N_n > j-1) = \sum_{j = 1}^\infty \sum_{i = j}^\infty q (i, n)= \sum_{j = 1}^\infty \sum_{i = j}^\infty ( F_\Lambda (i / n) - F_\Lambda ((i-1) / n))=\sum_{j = 0}^\infty \overline{F}_\Lambda (j / n).
\end{equation*}
Thus,
\begin{equation*}
\overline{C}_{0, n} = \sum_{j = 0}^\infty \overline{F}_\Lambda (j / n) / n.
\end{equation*}
Also,
\begin{equation*}
\overline{F}_{N_n}(i-1) = P ( N_n > i-1)= \sum_{k = i}^\infty q (k, n)= \sum_{k = i}^\infty ( F_\Lambda (k / n) - F_\Lambda ((k-1) / n))= \overline{F}_\Lambda ((i-1) / n).
\end{equation*}
Then,
\begin{equation*}
f_{Ne}(i) = \frac{\overline{F}_\Lambda ((i-1) / n)} {\sum_{j = 0}^\infty \overline{F}_\Lambda (j / n)},
\quad i\ge1.
\end{equation*} 
\qed

\subsection{First approximation method}

Our first proposal of approximation to $\psi(u)$ is presented as a corollary of
Theorem~\ref{principal-myb-theorem}.
Note that $ \overline{C}_{0, n} = \sum_{j = 0}^\infty \overline{F}_\Lambda (j / n) / n $
is an upper sum of the integral of $ \overline{F}_\Lambda $.
Thus,  $ \overline{C}_{0, n} \to E (\Lambda) $ as $ n \rightarrow \infty $.
For the approximation methods we propose, we will take $ \overline{C}_{0, n} = E (\Lambda) $,
for any value of $n$.
\begin{corollary} \label{coro-uno-cap5}
Suppose a Gerber-Dickson model with $ \mbox{MP}(F_\Lambda)$ claims  is given. For large $n$,
\begin{equation}
\label{approach-one-d}
\psi(u) \approx \sum_{k = 0}^\infty \overline{C}_{k, n} \cdot \mbox{nb}(u, 1 / (1 + n)) (k),
\end{equation}
where
\begin{eqnarray}
\label{C0A1} \overline{C}_{0, n} &=& E (\Lambda), \\
\label{CkA1} \overline{C}_{k, n} &=& E (\Lambda) \left [\sum_{i = 1}^{k} f_{Ne}(i)\,\overline{C}_{k-i, n} + \overline{F}_{Ne}(k) \right],\quad k\ge1, \\
\label{fNeA1} f_{Ne}(i) &=& \frac{\overline{F}_\Lambda ((i-1) / n)} {\sum_{j = 0}^\infty \overline{F}_\Lambda (j / n)},\quad i\ge1.
\end{eqnarray}
\end{corollary}
For the examples shown in the next section, we have numerically found that the sum
in~(\ref{approach-one-d}) quickly converge to its value.
%This occurs since the summands are the product of the probabilities of the negative binomial distribution
%and the coefficients $\overline{C}_{k, n} $, the later being proportional to tail
%probabilities~(\ref{560-como-proba-de-cola}).
This will allow us to truncate the infinite sum without much loss of accuracy.\\

For example, suppose claims have a $ \mbox{MP}(F_\Lambda ) $ distribution, where $ \Lambda \sim \mbox{Exp}(\beta) $. In this case, claims have $ \mbox{Geo}(\beta / (1+ \beta)) $ distribution and by (\ref{ruina-geo}),
\begin{equation*}
\psi (u) = \left(\frac{1 / (1+ \beta)} {\beta / (1+ \beta)} \right)^{u + 1} = \frac{1}{\beta^{u + 1}}.
\end{equation*}
We will check that our approximation (\ref{approach-one-d}) converges to this solution
as $ n\rightarrow \infty $. First, the following is easily calculated:
$f_{Ne}(i) = e^{- i \beta / n}(e^{\beta / n} -1)$ and $\overline{F}_{Ne}(k) = e^{-\beta k / n}$.
After some more calculations, one can obtain
\begin{equation}
\label{Cjn-d}
\overline {C}_{k, n} = \frac{1} {\beta} \left [\frac{1} {\beta}(1-e^{- \beta / n}) + e^{- \beta / n} \right]^k.
\end{equation}
Substituting (\ref{Cjn-d}) into (\ref{approach-one-d}) and simplifying,
\begin{displaymath}
\psi_n (u) 
= \frac{1} {\beta} \, \left(1-n \, (1-e^{- \beta / n}) / \beta + n \, (1-e^{- \beta / n}) \right)^{- u}
\ \to\ 1 / \beta^{u + 1}\quad\mbox{as}\quad n\to\infty.
\end{displaymath}
%t can be verified that $ n \, (1-e^{- \beta / n}) \rightarrow \beta $ as $n\to\infty$.
%Therefore,
%\begin{equation*}
%\lim_{n \rightarrow \infty} \psi_n (u)
%= \frac{1} {\beta} \, \left(1-1 + \beta \right)^{- u} 
%= 1 / \beta^{u + 1}.
%\end{equation*}
%which is the known ruin probability when claims are geometric distributed.

\subsection{Second approximation method}

Our second method to approximate the ruin probability is a direct
application of the Law of Large Numbers.

%, which is stated in the following corollary, is based on the expected value that appears in (\ref{560}), applying the Large Numbers Law (LNL) \cite[ p . 260] {1968feller}. (creo que no es necesario dar
% la referencia)

\begin{corollary} \label{coro-dos-cap5}
Suppose a Gerber-Dickson model with $ \mbox{MP}(F_\Lambda)$ claims  is given.
Let $ z_1, \ldots, z_m $ be a random sample of a $ \mbox{NegBin}(u, 1 / (1 + n)) $ distribution.
For large $ n $ and $ m $,
\begin{equation}
\label{two-d-approach}
\psi (u) \approx \frac{1} {m} \sum_{i = 1}^m \overline{C}_{z_i, n},
\end{equation}
where $ \{\overline{C}_{k, n} \}_{k = 0}^{\infty} $ is given by (\ref{C0A1}), (\ref{CkA1}) and (\ref{fNeA1}).
\end{corollary}

%-----------------------------------------------------------------

\section{Numerical examples}

In this section we apply the proposed approximation methods in the
case when the mixing distribution is Erlang, Pareto and Lognormal.
The results obtained show that the approximated ruin probabilities are
extremely close to the exact probabilities. The later were calculated
recursively using formulas~(\ref{psicerod}) and~(\ref{psid}), or by
numerical integration.
%The Erlang, Pareto and Lognormal distributions were used to model the mixing distribution of the MP claims.
%To model the amount of the claims, the exponential, Erlang , uniform, Pareto , lognormal and geometric truncated to zero distributions were used as the mixture distribution . 
In all cases the approximations were calculated for $ u = 0,1,2, \ldots, 10 $ and
using the software {\tt R}. For the first proposed approximation method, $ n = 500 $ was used and for the second method, $ m = 1000 $ values ​​were generated from a $ \mbox{NegBin}(u, 1 / (n + 1)) $ distribution and again $ n = $ 500. The
sum~(\ref{approach-one-d})  was truncated up to
\begin{equation}
k^* = \max \{x> un : \mbox{nb}(u, 1 / (1 + n)) (x)> 0.00001 \}.
\end{equation}
%The implementation of the methods was done using the R \cite{2008R} software.

\subsection*{Erlang  distribution}
\begin{figure}[H]
\centering
\input{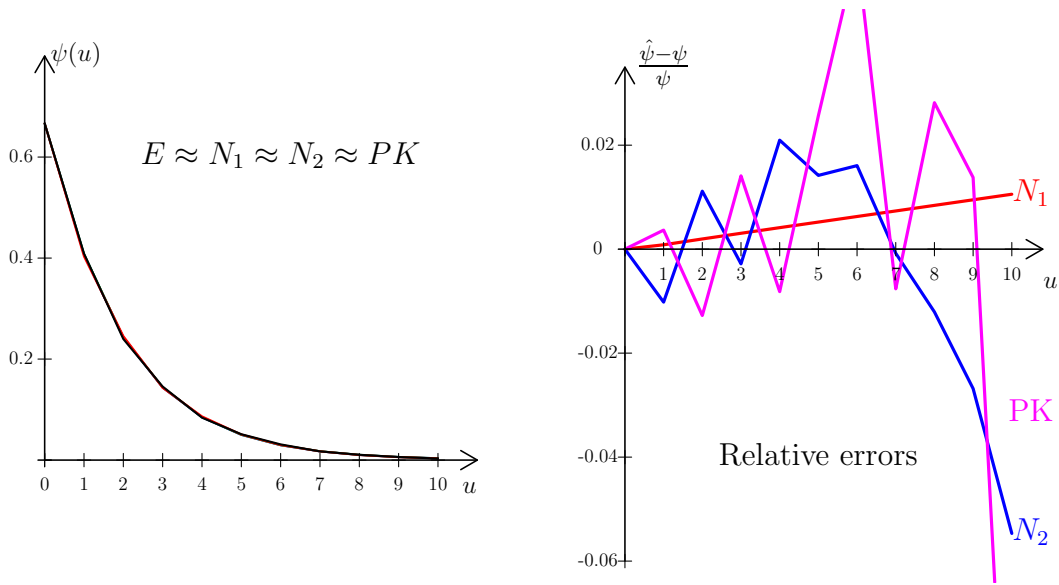}
\caption{Approximation when claims are $ \mbox{MP}(\Lambda) $ and $ \Lambda \sim \mbox{Erl}(2,3) $.} \label{dg}
\end{figure}

In this example we assume claims have a $ \mbox{MP}(F_\Lambda) $
distribution with $ \Lambda \sim \mbox{Erlang}(2,3) $.
In this case $ E (\Lambda) = 2/3 $. Table \ref{table-lg} below shows the results of the approximations. Columns $ E $, $ N_1 $ and $ N_2 $ show for each value of $ u $, the exact value of $\psi(u)$, the approximation with the first method and the approximation with the second method, respectively. Relative errors $ (\hat{\psi} - \psi) / \psi $ are also shown. The left-hand side plot of Figure~\ref{dg} shows
the values ​​of $ u $ against $ E $, $ N_1 $ and $ N_2 $. The right-hand side plot
shows the values ​​of $ u $ against the relative errors.

\subsection*{Pareto  distribution}
\begin{figure}[H]
\centering
\input{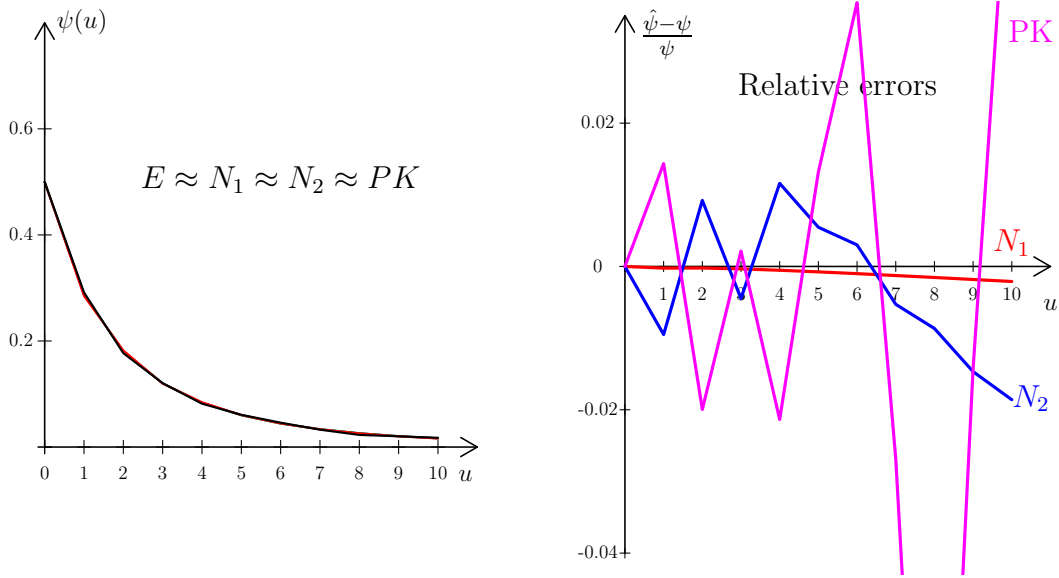}
\caption{Approximation when claims are $ \mbox{MP}(\Lambda) $ and $ \Lambda \sim \mbox{Pareto}(3,1) $.} \label{dp}
\end{figure}

In this example claims have a $ \mbox{MP}(F_\Lambda) $ distribution
with $ \Lambda \sim \mbox{Pareto}(3,1) $. For this distribution, $ E (\Lambda) = 1/2 $. Table \ref{table-lp} shows the approximations results in the same terms as in Table \ref{table-lg}. Figure \ref{dp} shows the results graphically.

\subsection*{Lognormal  distribution}
\begin{figure}[H]
\centering
\input{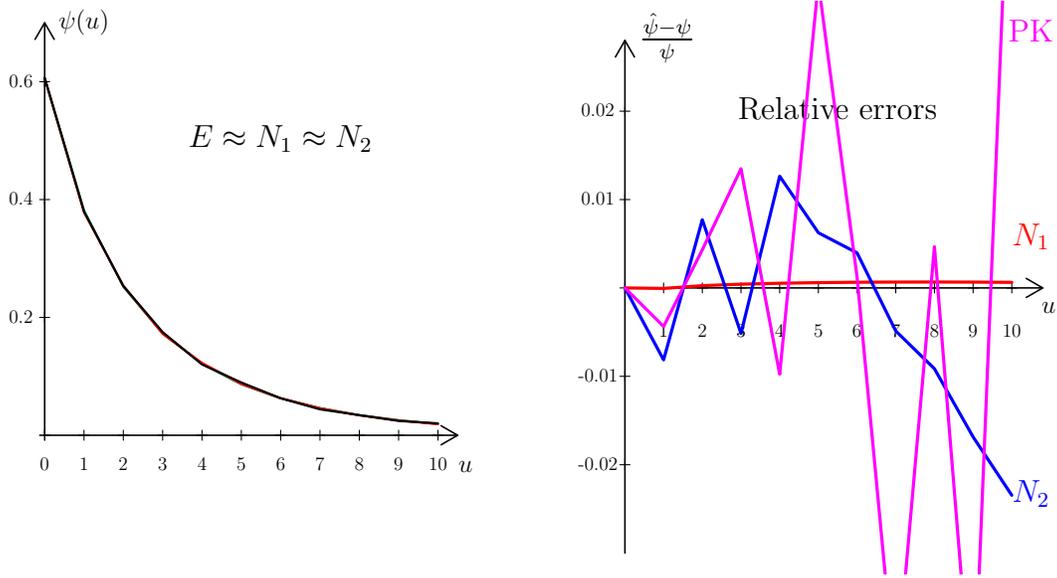}
\caption{Approximation when claims are $ \mbox{MP}(\Lambda) $ and $ \Lambda \sim \mbox{Lognormal}(- 1,1) $.} \label{dl}
\end{figure}

In this example we suppose claims have a $ \mbox{MP}(F_\Lambda) $ distribution with $ \Lambda \sim \mbox{Lognormal}(- 1,1) $. For this distribution $ E (\Lambda) = e^{- 1/2} $. Table \ref{table-ll} shows the approximations results and Figure \ref{dl} shows the related graphics.\\

As can be seen from the tables and graphs shown, the two approximating methods
yield ruin probabilities close to the exact probabilities for the examples
considered. 
\newpage

\begin{table}[H]  \large \centering
\caption{Ruin probability approximation for 
$ \mbox{MP}(F_\Lambda) $ claims with $ \Lambda \sim \mbox{Erlang}(2,3) $.}

\label{table-lg}

\begin{tabular} {@ {} cccccccc @ {}}

\toprule                                                                                               
$ u $ & $ E $ & $ N_1 $ & $ \frac{\hat{\psi} - \psi} {\psi} $ & $ N_2 $ & $ \frac{\hat{\psi} - \psi} {\psi} $ & $ PK $ & $ \frac{\hat{\psi} - \psi} {\psi} $\\ \midrule
0	&	0.66667	&	0.66667	&	0.00000	&	0.66667	&	0.00000	&	0.66667	&	0.00000	\\
1	&	0.40741	&	0.40775	&	0.00084	&	0.40326	&	-0.01019	&	0.4089	&	0.00366	\\
2	&	0.24280	&	0.24328	&	0.00196	&	0.24551	&	0.01115	&	0.2397	&	-0.01276	\\
3	&	0.14358	&	0.14401	&	0.00306	&	0.14317	&	-0.00282	&	0.1456	&	0.01410	\\
4	&	0.08469	&	0.08504	&	0.00414	&	0.08647	&	0.02096	&	0.084	&	-0.00818	\\
5	&	0.04992	&	0.05018	&	0.00521	&	0.05063	&	0.01419	&	0.0512	&	0.02566	\\
6	&	0.02942	&	0.02960	&	0.00628	&	0.02989	&	0.01607	&	0.0311	&	0.05726	\\
7	&	0.01733	&	0.01746	&	0.00735	&	0.01732	&	-0.00079	&	0.0172	&	-0.00763	\\
8	&	0.01021	&	0.01030	&	0.00842	&	0.01009	&	-0.01208	&	0.0105	&	0.02818	\\
9	&	0.00602	&	0.00607	&	0.00949	&	0.00586	&	-0.02682	&	0.0061	&	0.01379	\\
10	&	0.00355	&	0.00358	&	0.01056	&	0.00335	&	-0.05468	&	0.0031	&	-0.12559	\\ \bottomrule
\end{tabular}
\end{table}                    

\begin{table}[H]  \large \centering
\caption{Ruin probability approximation for $ \mbox{MP}(F_\Lambda) $ claims
with $ \Lambda \sim \mbox{Pareto}(3,1) $.}
\label{table-lp}
\begin{tabular} {@ {} cccccccc @ {}}
\toprule                                                                                               
$ u $ & $ E $ & $ N_1 $ & $ \frac{\hat{\psi} - \psi} {\psi} $ & $ N_2 $ & $ \frac{\hat{\psi} - \psi} {\psi} $ & $ PK $ & $ \frac{\hat{\psi} - \psi} {\psi} $\\ \midrule
0	&	0.50000	&	0.50000	&	0.00000	&	0.50000	&	0.00000	&	0.50000	&	0.00000	\\
1	&	0.28757	&	0.28751	&	-0.00023	&	0.28484	&	-0.00950	&	0.29170	&	0.01435	\\
2	&	0.18050	&	0.18046	&	-0.00022	&	0.18216	&	0.00921	&	0.17690	&	-0.01995	\\
3	&	0.12014	&	0.12010	&	-0.00034	&	0.11960	&	-0.00448	&	0.12040	&	0.00215	\\
4	&	0.08348	&	0.08344	&	-0.00053	&	0.08445	&	0.01159	&	0.08170	&	-0.02135	\\
5	&	0.06001	&	0.05996	&	-0.00076	&	0.06034	&	0.00547	&	0.06080	&	0.01317	\\
6	&	0.04437	&	0.04432	&	-0.00100	&	0.04450	&	0.00301	&	0.04600	&	0.03681	\\
7	&	0.03360	&	0.03356	&	-0.00127	&	0.03343	&	-0.00528	&	0.03270	&	-0.02686	\\
8	&	0.02599	&	0.02595	&	-0.00154	&	0.02577	&	-0.00865	&	0.02280	&	-0.12288	\\
9	&	0.02049	&	0.02045	&	-0.00181	&	0.02019	&	-0.01467	&	0.02020	&	-0.01419	\\
10	&	0.01643	&	0.01639	&	-0.00209	&	0.01612	&	-0.01858	&	0.01750	&	0.06527	\\

\bottomrule

\end{tabular}

\end{table}

\begin{table}[H]  \large \centering
\caption{Approximations for $ \mbox{MP}(F_\Lambda) $ claims
with $ \Lambda \sim \mbox{Lognormal}(-1,1) $.}
\label{table-ll}
\begin{tabular} {@ {} cccccccc @ {}}
\toprule                                                                                               
$ u $ & $ E $ & $ N_1 $ & $ \frac{\hat{\psi} - \psi} {\psi} $ & $ N_2 $ & $ \frac{\hat{\psi} - \psi} {\psi} $ & $ PK $ & $ \frac{\hat{\psi} - \psi} {\psi} $\\ \midrule
0	&	0.60653	&	0.60653	&	0.00000	&	0.60653	&	0.00000	&	0.60653	&	0.00000	\\
1	&	0.38126	&	0.38124	&	-0.00005	&	0.37816	&	-0.00813	&	0.37960	&	-0.00436	\\
2	&	0.25231	&	0.25238	&	0.00025	&	0.25426	&	0.00772	&	0.25340	&	0.00431	\\
3	&	0.17287	&	0.17294	&	0.00042	&	0.17198	&	-0.00515	&	0.17520	&	0.01349	\\
4	&	0.12128	&	0.12135	&	0.00053	&	0.12282	&	0.01264	&	0.12010	&	-0.00976	\\
5	&	0.08661	&	0.08666	&	0.00060	&	0.08715	&	0.00624	&	0.08960	&	0.03456	\\
6	&	0.06272	&	0.06276	&	0.00064	&	0.06297	&	0.00397	&	0.06280	&	0.00124	\\
7	&	0.04597	&	0.04600	&	0.00067	&	0.04574	&	-0.00487	&	0.04390	&	-0.04498	\\
8	&	0.03404	&	0.03406	&	0.00067	&	0.03373	&	-0.00914	&	0.03420	&	0.00466	\\
9	&	0.02545	&	0.02546	&	0.00066	&	0.02502	&	-0.01686	&	0.02420	&	-0.04902	\\
10	&	0.01919	&	0.01920	&	0.00063	&	0.01874	&	-0.02346	&	0.02030	&	0.05791	\\

 \bottomrule
\end{tabular}
\end{table}

\section{Conclusions}
\label{conclusions}

We have first provided a general formula for the ultimate ruin probability in the
Gerber-Dickson risk model when claims follow a negative binomial mixture (NBM)
distribution. The ruin probability is expressed as the expected value of a deterministic
sequence $\{C_k\}$, where index $k$ is the value of a negative binomial distribution. The
sequence is not given explicitly but can be calculated recursively. We then extended
the formula for claims with a mixed Poisson (MP) distribution. The extension was
possible due to the fact that MP distributions can be approximated by NBM distributions.
The formulas obtained yielded two immediate approximation methods. These were tested
using particular examples. The numerical results showed high accuracy when compared
to the exact ruin probabilities. The general results obtained in this work bring about some
other questions that we have set aside for further work: error bounds for our estimates,
detailed study of some other particular cases of the NBM and MP distributions, properties
and bounds for the sequence $\{C_k\}$, and the possible extension of the ruin probability
formula to more general claim distributions.

\end{document}